\documentclass[smallextended,envcountsame,envcountsect]{svjour3}

\smartqed
\usepackage{mathptmx}

\usepackage[utf8]{inputenc}
\usepackage{tikz}
\usepackage{amssymb,amsmath, graphicx,hyperref, xcolor,lineno}
\usepackage{newtxtext,newtxmath}
\usepackage{parskip}

\newcommand{\C}{\mathcal{C}}
\newcommand{\aw}{\textup{aw}}
\newcommand{\dist}{\textup{d}}
\newcommand{\diam}{\textup{diam}}
\newcommand{\Mod}[1]{\ (\mathrm{mod}\ #1)}

\title{Anti-van der Waerden Numbers of Graph Products of Cycles\thanks{Thank you to the University of Wisconsin-La Crosse's (UWL's) Dean's Distinguished Fellowship program that supported both authors.  Also, thanks to the UWL's Undergraduate Research and Creativity grant and the UWL Department of Mathematics and Statistics Bange/Wine Undergraduate Research Endowment that supported the first author.}
}

\author{Joe Miller \and Nathan Warnberg}

\institute{University of Wisconsin-La Crosse\\ \email{miller1374@uwlax.edu \and nwarnberg@uwlax.edu}}

\institute{J. Miller \and N. Warnberg \\
              University of Wisconsin-La Crosse \\
              \email{miller1374@uwlax.edu \and nwarnberg@uwlax.edu} }

\date{\today}

\begin{document}


\titlerunning
\authorrunning

\maketitle

\begin{abstract}

A \emph{$k$-term arithmetic progression ($k$-AP) in a graph $G$} is a list of vertices such that each consecutive pair of vertices is the same distance apart.  If $c$ is a coloring function of the vertices of $G$ and a $k$-AP in $G$ has each vertex colored distinctly, then that $k$-AP is a \emph{rainbow $k$-AP}. The \emph{anti-van der Waerden number of a graph $G$ with respect to $k$} is the least positive integer $r$ such that every surjective coloring with domain $V(G)$ and codomain $\{1,2,\dots,r\} = [r]$ is guaranteed to have a rainbow $k$-AP.  This paper focuses on $3$-APs and graph products with cycles.  Specifically, the anti-van der Waerden number with respect to $3$ is determined precisely for $P_m \square C_n$, $C_m\square C_n$ and $G\square C_{2n+1}$.


\end{abstract}

\newpage

\section{Introduction}\label{sec:intro}

The study van der Waerden numbers started with Bartel van der Waerden showing in $1927$ that given a fixed number of colors $r$, and a fixed integer $k$ there is some $N$ (a van der Waerden number) such that if $n \ge N$, then no matter how you color $[n] = \{1,2,\dots,n\}$ with $r$-colors, there will always be a monochromatic $k$-term arithmetic progression (see \cite{W27}).  Around this time, in $1917$, it is interesting to note that I. Schur proved that given $r$ colors, you can find an $N$ (a Schur number) such that if $n \ge N$, then no matter how you color $[n]$ there must be a monochromatic solution to $x+y = z$ (see \cite{S}).  In addition, in $1928$, F.P. Ramsey showed that (here graph theory language is used but was not in Ramsey's original formulation) given $r$ colors and some constant $k$ you can find an $N$ (a Ramsey number) such that if $n \ge N$, then no matter how you colors the edges of a complete graph $K_n$ you can always find a complete subgraph $K_k$ that is monochromatic (see \cite{R}).  

These types of problems that look for monochromatic sturctures have been categorized as Ramsey-type problems and each of them have a dual version.  For example, an anti-van der Waerden number is when given integers $n$ and $k$, find the minimum number of colors such that coloring $\{1,\dots,n\}$ ensures a rainbow $k$-term arithmetic progression.  It was not until 1973 when Erd\H{o}s, Simonovits, and S\'{o}s, in \cite{ESS}, started looking at the dual versions of these problems which are now well-studied (see \cite{FMO} for a survey).

Results on colorings and balanced colorings of $[n]$ that avoid rainbow arithmetic progressions have been studied in \cite{AF} and \cite{AM}.  Rainbow free colorings of $[n]$ and $\mathbb{Z}_n$ were studied in \cite{J} and \cite{DMS}.  Although Butler et al., in \cite{DMS}, consider arithmetic progressions of all lengths, many results on $3$-APs were produced.  In particular, the authors of \cite{DMS} determine $\aw(\mathbb{Z}_n,3)$, see Theorem \ref{cycles} with additional cycle notation.  Further, the authors of \cite{DMS} determined that $3\leq \aw(\mathbb{Z}_p,3)\leq 4$ for every prime number $p$ and that $\aw(\mathbb{Z}_n,3)$, which equals $\aw(C_n,3)$, can be determined by the prime factorization of $n$.  This result was then generalized by Young in \cite{finabgroup}. 

\begin{theorem}\label{cycles} {\cite{DMS}}
\begin{sloppy}
Let $n$ be a positive integer with prime decomposition $n=2^{e_0}p_1^{e_1}p_2^{e_2}\cdots p_s^{e_s}$ for $e_i\geq 0$, $i=0,\ldots,s$, where primes are ordered so that $\aw(\mathbb{Z}_{p_i},3)=3$ for $ 1 \leq i \leq \ell$ and $\aw(\mathbb{Z}_{p_i},3)=4$ for $\ell + 1 \leq i \leq s$. Then, 
	\[
	\aw(\mathbb{Z}_n,3) = 		\aw(C_n,3)=\left\{\begin{array}{ll}
				2 +\sum\limits_{j=1}^\ell e_j + \sum\limits_{j=\ell+1}^s 2e_j & \mbox{if $n$ is odd,} \\
				3 +\sum\limits_{j=1}^\ell e_j + \sum\limits_{j=\ell+1}^s 2e_j & \mbox{if $n$ is even.}
			\end{array}\right.
		\]	
		\end{sloppy}
	\end{theorem}
	
As mentioned, Butler et al. also studied arithmetic progressions on $[n]$ and obtained bounds on $\aw([n],3)$ and conjectured the exact value that was later proven in \cite{BSY}.  This result on $[n]$ is presented as Theorem \ref{paths} and includes path notation.

	\begin{theorem}\label{paths}\cite{BSY}
If $n \ge 3$ and $7\cdot 3^{m-2} +1 \le n \le 21\cdot 3^{m-2}$, then
	\[\aw([n],3) = \aw(P_n,3) = \left\{\begin{array}{ll} m+2 & \text{ if $n = 3^m$,}\\
											m+3 & \text{ otherwise.}\end{array}\right. \]
									
\end{theorem}
	
It is also interesting to note that $3$-APs in $[n]$ or $\mathbb{Z}_n$ satisfy the equation $x_1 + x_2 = 2x_3$. Thus, rainbow numbers for other linear equations have also been considered (see \cite{BKKTTY}, \cite{FGRWW}, \cite{RFC} and \cite{LM}).

Studying the anti-van der Waerden number of graphs is a natural extension of determining anti-van der Waerden number of $[n] = \{1,2,\dots,n\}$, which behave like paths, and $\mathbb{Z}_n$, which behave like cycles.  In particular, the set of arithmetic progressions on $[n]$ is isomorphic to the set of arithmetic progressions on $P_n$ and the set of arithmetic progressions on $\mathbb{Z}_n$ is isomorphic to the set of arithmetic progressions on $C_n$. This relationship was first introduced and explored in \cite{SWY} where the anti-van der Waerden number was bounded by the radius and diameter of a graph, the anti-van der Waerden number of trees and hypercubes were investigated and an upperbound of $4$ was conjectured for the anti-van der Waerden number of graph products. Then, in \cite{RSW}, the authors confirmed the upper bound of $4$ for any graph product (see Theorem \ref{thm:rsw}).  This paper continues in this vein.

\begin{theorem}\cite{RSW}\label{thm:rsw}
    If $G$ and $H$ are connected graphs and $|G|,|H| \ge 2$, then \[\aw(G\square H,3) \le 4.\]
\end{theorem}

Something that makes anti-van der Waerden numbers challenging is that it is not a subgraph monotone parameter. A particular example, \[4 = \aw([9],3) = \aw(P_9,3) < \aw(P_8,3) = \aw([8],3) = 5,\] even though $P_8$ is a subgraph of $P_9$, and a general statement, \[\aw(C_n,3) = \aw(\mathbb{Z}_n,3) \le \aw([n],3) = \aw(P_n,3),\] were both given, without the graph theory interpretation, in \cite{DMS}.  One tool that does allow a kind of monotonicity when studying the anti-van der Waerden numbers of graphs is when a subgraph is isometric, that is, the subgraph preserves distances.  This insight was used extensively in \cite{RSW} to get an upper bound on the anti-van der Waerden number of graph products and will also be leveraged in this paper.  First, some definitions and background inspired by \cite{DMS} and used in \cite{SWY} and \cite{RSW} are provided.

Graphs in this paper are undirected so edge $\{u,v\}$ will be shortened to $uv\in E(G)$.  If $uv\in E(G)$, then $u$ and $v$ are \emph{neighbors} of each other. The \emph{distance} between vertex $u$ and $v$ in graph $G$ is denoted $\dist_G(u, v)$, or just $\dist(u, v)$ when context is clear, and is the smallest length of any $u-v$ path in $G$.  A $u-v$ path of length $d(u,v)$ is called a $u-v$ \emph{geodesic}.

A \emph{$k$-term arithmetic progression in graph $G$} ($k$-AP) is a set of vertices $\{v_1,v_2,\dots,v_k\}$ such that $d(v_i,v_{i+1}) = d$ for all $1\le i \le k-1$.  A $k$-term arithmetic is \emph{degenerate} if $v_i = v_j$ for any $i\neq j$.  Note that technically, since a $k$-AP is a set, the order of the elements does not matter.  However, oftentimes $k$-APs will be presented in the order that provides the most intuition.  

An \emph{exact $r$-coloring of a graph $G$} is a surjective function $c:V(G) \to [r]$.  A set of vertices $S$ is \emph{rainbow} under coloring $c$ if for every $v_i,v_j\in V(G)$, $c(v_i) \neq c(v_j)$ when $i\neq j$.  Given a set $S\subset V(G)$, define $c(S) = \{ c(s) \, | \, s\in S\}$.

The \emph{anti-van der Waerden number of graph $G$ with respect to $k$}, denoted $\aw(G,k)$, is the least positive number $r$ such that every exact $r$-coloring of $G$ contains a rainbow $k$-term arithmetic progression.  If $|V(G)| = n$ and no coloring of the vertices yields a rainbow $k$-AP, then $\aw(G,k) = n+1$.

Graph $G'$ is a \emph{subgraph} of $G$ if $V(G') \subseteq V(G)$
or $E(G') \subseteq E(G)$ (or both). A subgraph $G'$ of $G$ is an \emph{induced subgraph} if whenever $u$ and $v$ are vertices of $G'$ and $uv$ is an edge of $G$, then $uv$ is an edge of $G'$.  If $S$ is a nonempty set of vertices of $G$, then the \emph{subgraph of $G$ induced by $S$} is the induced subgraph with vertex set $S$ and is denoted $G[S]$.  An \emph{isometric subgraph} $G'$ of $G$ is a subgraph such that $\dist_{G'}(u,v) = \dist_G(u,v)$ for all $u,v \in V(G')$. 

If $G = (V,E)$ and $H = (V', E')$ then the \emph{Cartesian
product}, written $G\square H$, has vertex set $\{(x, y) : x \in V \text{ and } y \in V' \}$ and $(x, y)$ and
$(x', y')$ are adjacent in $G \square H$ if either $x = x'$ and $yy' \in E'$ or $y = y'$ and $xx' \in E$. 

This paper will use the convention that if \[V(G) = \{u_1,\ldots, u_{n_1}\} \quad \text{and} \quad V(H) = \{w_1,\ldots,w_{n_2}\},\] then $V(G\square H) = \{v_{1,1},\ldots, v_{n_1,n_2}\}$ where $v_{i,j}$ corresponds to the vertices $u_i \in V(G)$ and $w_j \in V(H)$.

Also, if $1\leq i \leq n_2$, then $G_i$ denotes the $i$th labeled copy of $G$ in $G \square H$. Likewise, if $1 \leq j \leq n_1$, then $H_j$ denotes the $j$th labeled copy of $H$ in $G \square H$.  In other words, $G_i$ is the induced subgraph $G_i = G\square H[\{v_{1,i},\dots, v_{n_2,i}\}]$, and $H_j$ is the induced subgraph $H_j = G\square H[\{v_{j,1}, \dots, v_{j,n_1}\}]$.  Notice that the $i$ subscript in $G_i$ corresponds to the $i$th vertex of $H$ and the $j$ in the subscript in $H_j$ corresponds to the $j$th vertex of $G$.  See Example \ref{ex:cartprod} below.

\begin{example}\label{ex:cartprod}

  Consider the graph $P_3\square C_5$ where $V(P_3) = \{u_1,u_2,u_3\}$ and $V(C_5) = \{w_1,w_2,w_3,w_4,w_4\}$.  Let $G=P_3$ and $H = C_5$ as in the definition.  Now, $G_4$ is a subgraph of $P_3\square C_5$ that is isomorphic to $P_3$ and corresponds to vertex $w_4$ of $C_5$.  Similarly, $H_2$ is a subgraph of $P_3\square C_5$ that is isomorphic to $C_5$ and corresponds to vertex $u_2$ of $P_3$.  See Figure \ref{fig:cartprodex} below.
  
  \begin{figure}[ht!]
      \centering
      \begin{tikzpicture}[scale = 0.7]
        \node[draw,circle] (11) at (.75,0) {$v_{1,1}$};
	    \node[draw,circle] (12) at (3,0) {$v_{1,2}$};
	    \node[draw,circle] (13) at (3.75,2.25) {$v_{1,3}$};
	    \node[draw,circle, line width=0.75mm] (14) at (1.9875,3.75) {$v_{1,4}$};
	    \node[draw,circle] (15) at (0,2.25) {$v_{1,5}$};
	    
	    \draw[thick]  (11) to node [auto] {} (12);
      	\draw[thick]  (12) to node [auto] {} (13);
      	\draw[thick]  (13) to node [auto] {} (14);
      	\draw[thick]  (14) to node [auto] {} (15);
      	\draw[thick]  (15) to node [auto] {} (11);
      	
      	\node[draw,circle, thick, dashed] (21) at (5.75,1) {$v_{2,1}$};
	    \node[draw,circle, thick, dashed] (22) at (8,1) {$v_{2,2}$};
	    \node[draw,circle, thick, dashed] (23) at (8.75,3.25) {$v_{2,3}$};
	    \node[draw,circle, line width=0.75mm, dashed] (24) at (6.9875,4.75) {$v_{2,4}$};
	    \node[draw,circle, thick, dashed] (25) at (5,3.25) {$v_{2,5}$};
	    
	    \draw[thick, dashed]  (21) to node [auto] {} (22);
      	\draw[thick, dashed]  (22) to node [auto] {} (23);
      	\draw[thick, dashed]  (23) to node [auto] {} (24);
      	\draw[thick, dashed]  (24) to node [auto] {} (25);
      	\draw[thick, dashed]  (25) to node [auto] {} (21);
      	
      	\node[draw,circle] (31) at (10.75,2) {$v_{3,1}$};
	    \node[draw,circle] (32) at (13,2) {$v_{3,2}$};
	    \node[draw,circle] (33) at (13.75,4.25) {$v_{3,3}$};
	    \node[draw,circle, line width=0.75mm] (34) at (11.9875,5.75) {$v_{3,4}$};
	    \node[draw,circle] (35) at (10,4.25) {$v_{3,5}$};
	    
	    \draw[thick]  (31) to node [auto] {} (32);
      	\draw[thick]  (32) to node [auto] {} (33);
      	\draw[thick]  (33) to node [auto] {} (34);
      	\draw[thick]  (34) to node [auto] {} (35);
      	\draw[thick]  (35) to node [auto] {} (31);
      	
      	\draw[thick]  (11) to node [auto] {} (21);
      	\draw[thick]  (12) to node [auto] {} (22);
      	\draw[thick]  (13) to node [auto] {} (23);
      	\draw[thick, line width=0.75mm]  (14) to node [auto] {} (24);
      	\draw[thick]  (15) to node [auto] {} (25);
      	
      	\draw[thick]  (31) to node [auto] {} (21);
      	\draw[thick]  (32) to node [auto] {} (22);
      	\draw[thick]  (33) to node [auto] {} (23);
      	\draw[thick, line width=0.75mm]  (34) to node [auto] {} (24);
      	\draw[thick]  (35) to node [auto] {} (25);
      	
      	\node[draw,circle,fill=white] (11) at (.75,0) {$v_{1,1}$};
	    \node[draw,circle,fill=white] (12) at (3,0) {$v_{1,2}$};
	    \node[draw,circle,fill=white] (13) at (3.75,2.25) {$v_{1,3}$};
	    \node[draw,circle, line width=0.75mm,fill=white] (14) at (1.9875,3.75) {$v_{1,4}$};
	    \node[draw,circle,fill=white] (15) at (0,2.25) {$v_{1,5}$};
      	
      	\node[draw,circle, thick, dashed,fill=white] (21) at (5.75,1) {$v_{2,1}$};
	    \node[draw,circle, thick, dashed,fill=white] (22) at (8,1) {$v_{2,2}$};
	    \node[draw,circle, thick, dashed,fill=white] (23) at (8.75,3.25) {$v_{2,3}$};
	    \node[draw,circle, line width=0.75mm, dashed,fill=white] (24) at (6.9875,4.75) {$v_{2,4}$};
	    \node[draw,circle, thick, dashed,fill=white] (25) at (5,3.25) {$v_{2,5}$};
      	
      	\node[draw,circle,fill=white] (31) at (10.75,2) {$v_{3,1}$};
	    \node[draw,circle,fill=white] (32) at (13,2) {$v_{3,2}$};
	    \node[draw,circle,fill=white] (33) at (13.75,4.25) {$v_{3,3}$};
	    \node[draw,circle, line width=0.75mm,fill=white] (34) at (11.9875,5.75) {$v_{3,4}$};
	    \node[draw,circle,fill=white] (35) at (10,4.25) {$v_{3,5}$};
      	
      \end{tikzpicture}
     
      \caption{Image for Example \ref{ex:cartprod}: The subgraph $G_4$ has been bolded and $H_2$ has been dashed.} 
      
      \label{fig:cartprodex}
  \end{figure}
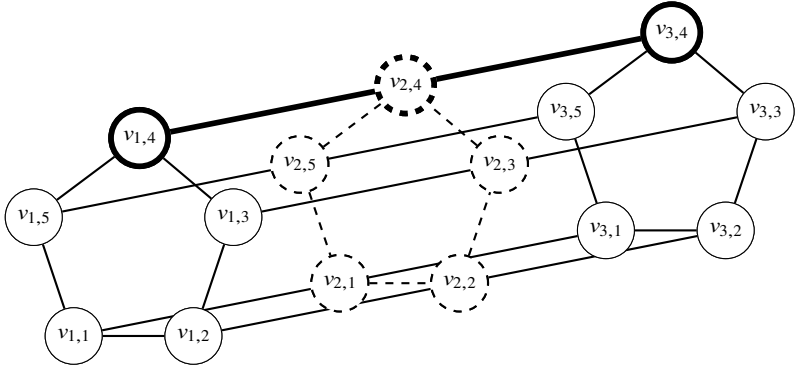
  
\end{example}

  The paper continues with Section \ref{sec:2} recapping and expanding many fundamental results from \cite{RSW}.  Section \ref{sec:pmcn} establishes $\aw(P_m\square C_n,3)$ for all $m$ and $n$.  Section \ref{sec:cmcn} is an investigation of $\aw(G\square C_n,3)$.  In particular, $\aw(C_m\square C_n,3)$ is determined for all $m$ and $n$.  Further, Section \ref{sec:cmcn} determines $\aw(G\square C_n,3)$ for any $G$ when $n$ is odd.  Finally, Section \ref{sec:future} provides the reader with some conjectures and open questions.

\section{Background and Fundamental Tools}\label{sec:2}

Distance preservation in subgraphs can be leveraged to guarantee the existence of rainbow $3$-APs. Thus, this section starts with some basic distance and isometry results.

\begin{proposition}\label{prop:dist}
If $v_{i,j},v_{h,k} \in V(G \square H)$, then 
\[\dist_{G\square H}(v_{i,j},v_{h,k}) = \dist_G(u_i,u_h) + \dist_H(w_j,w_k).\]
\end{proposition}

\begin{proof}
Note that $\dist_{G\square H}(v_{i,j},v_{h,k}) \le \dist_G(u_i,u_h) + \dist_H(w_j,w_k)$ because a path of length $\dist_G(u_i,u_h) + \dist_H(w_j,w_k)$ can be constructed using a $u_i-u_h$ geodesic in $G$ and combining it with a $w_j-w_k$ geodesic in $H$. 
\par
To show the other inequality, let $P$ be a $v_{i,j} - v_{h,k}$ geodesic, say \[P = \{v_{i,j} = x_1, x_2, \ldots, x_y = v_{h,k}\}.\] Note that for every edge $v_{j_1,j_2}v_{\beta_1,\beta_2} \in E(P)$, either $j_1 = \beta_1$ and $w_{j_2}w_{\beta_2} \in E(H)$, or $j_2 = \beta_2$ and $u_{j_1}u_{\beta_1} \in E(G)$. Then, $x_{\ell}x_{\ell+1}$ must correspond either to an edge from a $u_i-u_j$ walk or from a $w_h-w_k$ walk and $P$ must correspond to a walk in $G$ and also a walk in $H$. In other words, the length of $P$ is the sum of the length of the corresponding walks in $G$ and $H$. Thus, the length of $P$ is at least the sum of the lengths of a $u_i-u_h$ geodesic in $G$ and a $w_j-w_k$ geodesic in $H$. So, \[\dist_G(u_i,u_h) + \dist_H(w_j,w_k) \le \dist_{G\square H}(v_{i,j},v_{h,k}).\]
\qed\end{proof}

\begin{corollary}\label{cor:isosubprod}
    If $G'$ is an isometric subgraph of $G$ and $H'$ is an isometric subgraph of $H$, then $G'\square H'$  is an isometric subgraph of $G\square H$. 
\end{corollary}

\begin{proof}
    Let $V(G) = \{u_1,\ldots, u_{n_1}\}$ and $V(H) = \{w_1,\ldots,w_{n_2}\}$. Then let $v_{i,j},v_{h,k} \in V(G'\square H')$. Observe,
    \begin{align*}
        \dist_{G'\square H'}(v_{i,j},v_{h,k}) & = \dist_{G'}(u_i,u_h) + \dist_{H'}(w_j,w_k) \\
        & = \dist_{G}(u_i,u_h) + \dist_{H}(w_j,w_k) \\
        & = \dist_{G\square H}(v_{i,j},v_{h,k}).
    \end{align*}
    
\null\qed\end{proof}

Lemma \ref{isometricpathorC3} is powerful since it guarantees isometric subgraphs.  Isometric subgraphs are important when investigating anti-van der Waerden numbers because distance preservation implies $k$-AP preservation.

\begin{lemma}\cite{RSW}\label{isometricpathorC3}
    If $G$ is a connected graph on at least three vertices with an exact $r$-coloring $c$ where $r \ge 3$, then there exists a subgraph $G'$ in $G$ with at least three colors where $G'$ is either an isometric path or $G' = C_3$.
\end{lemma}

Theorem \ref{PmxPn} is used when isometric $P_m\square P_n$ subgraphs are found within $G\square H$.

\begin{theorem}\cite{RSW}\label{PmxPn}
    For $m,n \geq 2$, \[\aw(P_m \square P_n, 3) = \begin{cases}
    3 & \text{if $m = 2$ and $n$ is even, or $m = 3$ and $n$ is odd,} \\
    4 & \text{otherwise.}
    \end{cases}\]
\end{theorem}

Lemma \ref{|c(V(Gi U Gj))|<3} helps restrict the number of colors each copy of $G$ or $H$ can have within $G\square H$.

\begin{lemma}\cite{RSW}\label{|c(V(Gi U Gj))|<3}
    Assume $G$ and $H$ are connected with $|V(H)| \geq 3$. Suppose $c$ is an exact, rainbow-free $r$-coloring of $G\square H$, such that $r \geq 3$ and $|c(V(G_i))| \leq 2$ for $1 \leq i \leq n$. If $w_iw_j \in E(H)$, then $|c(V(G_i) \cup V (G_j))| \leq 2$.
\end{lemma}

To prove Lemmas \ref{lem:p2codd} and \ref{lem:p2ceven} requires the use of Lemma \ref{|c(Hi)|<3}.  

\begin{lemma}\label{|c(Hi)|<3}
    If $G$ and $H$ are connected, $|G|,|H| \ge 2$ and $c$ is an exact $r$-coloring of $G\square H$, $3\le r$, that avoids rainbow $3$-APs, then $|c(V(G_i))| \leq 2$ for $1 \leq i \leq |H|$.
\end{lemma}

\begin{proof} 
If $|G| = 2$ the result is immediate, so let $3\le |G|$.
For the sake of contradiction, assume $red,blue,green \in |c(V(G_i))|$ for some $1\le i \le |H|$. By Lemma \ref{isometricpathorC3}, there must exist an isometric path or a $C_3$ in $G_i$ containing $red$, $blue$, and $green$. If there is such a $C_3$, then there is a rainbow $3$-AP which is a contradiction. So, assume $P_\ell$ is a shortest isometric path in $G_i$ containing $red$, $blue$, and $green$, for some positive integer $3 \le \ell$. 
\begin{description}
\item[Case 1.] $\ell$ is odd.\\
Without loss of generality, suppose the two leaves of $P_\ell$ are colored $red$ and $blue$.  Since $P_\ell$ is shortest the rest of the vertices are colored $green$. Since $\ell$ is odd there exists a $green$ vertex equidistant from the $red$ and $blue$ vertices which creates a rainbow $3$-AP, a contradiction.
\item[Case 2.] $\ell$ is even.\\
Let $u_i \in V(H)$ be the vertex that corresponds to $G_i$ and note that $u_i$ has a neighbor since $H$ is connected. Let $P_2$ be a path on two vertices in $H$ containing $u_i$ and $\rho$ be the isometric subgraph in $G$ that corresponds to $P_\ell$. Thus, the subgraph $P_2 \square \rho$ of $G\square H$ is isometric and, by Theorem \ref{PmxPn}, contains a rainbow $3$-AP, a contradiction.
\end{description}
All cases give a contradiction, thus $|c(V(G_i))| \leq 2$.
\qed\end{proof}

Corollary \ref{cor:neighborcopies} is a strengthening of Lemma \ref{|c(V(Gi U Gj))|<3} and follows from Lemmas \ref{|c(V(Gi U Gj))|<3} and \ref{|c(Hi)|<3}.  It is used to help analyze $\aw(P_m\square C_{2k+1})$.

\begin{corollary}\label{cor:neighborcopies}
    If $G$ and $H$ are connected graphs, $|G| \ge 2$, $|H|\ge 3$, $c$ is an exact, rainbow-free $r$-coloring of $G\square H$ with $r\ge 3$, and $v_iv_j \in E(H)$, then \[|c(V(G_i)\cup V(G_j))| \leq 2.\]
\end{corollary}

\begin{lemma}\label{c(H_i)/c(H_j) < 2}\cite{SWY}
    Let $G$ be a connected graph on $m$ vertices and $H$ be a connected graph on $n$ vertices. Let $c$ be an exact $r$-coloring of $G\square H$ with no rainbow $3$-APs. If $G_1,G_2, \ldots,G_n$ are the labeled copies of $G$ in $G\square H$, then $|c(V(G_j)) \setminus c(V(G_i))| \leq 1$ for all $1 \leq i, j \leq n$.
\end{lemma}

\begin{proposition}\label{prop:everycopy}
    If $G$ and $H$ are connected graphs, $|G| \ge 2$, $|H|\ge 3$, $c$ is an exact, rainbow-free $r$-coloring of $G\square H$ with $r\ge 3$, then there is a color in $c(G\square H)$ that appears in every copy of $G$.
\end{proposition}

\begin{proof}
Suppose $c(G\square H) = \{c_1,\ldots,c_r\}$. First, for the sake of contradiction, assume $|c(V(G_i))|=1$ for every $1 \leq i \leq |H|$. Then define a coloring $c': V(H) \to c(G\square H)$ such that $c'(w_i) \in c(V(G_i))$. Then Lemma \ref{isometricpathorC3} implies that there is either an isometric path or $C_3$ in $H$ with $3$ colors. If there is an isometric $C_3$, say $(w_1,w_2,w_3)$, then $\{v_{1,1},v_{1,2},v_{1,3}\}$ is a rainbow $3$-AP in $G\square H$ with respect to $c$, a contradiction. So, there must be an isometric path in $H$ with $3$ colors. Suppose $P = (w_1,\ldots,w_n)$ is a shortest such path. Without loss of generality, $c(w_1)=c_2$, $c(w_n)=c_3$ and $c(w_i)=c_1$ for all $1 < i < n$. Then there exists $u_1,u_2\in V(G)$ such that $u_1u_2 \in E(G)$. Thus, $\{v_{1,1},v_{1,n},v_{2,2}\}$ is a rainbow $3$-AP in $G\square H$ with respect to $c$, a contradiction.

Thus, there exists some $G_i$ such that $|c(V(G_i))|\geq 2$, without loss of generality, say $c_1,c_2 \in c(V(G_i))$. Then Lemma \ref{|c(Hi)|<3} implies $c(V(G_i))=\{c_1,c_2\}$. Note that $c_3 \in c(V(G_j))$ for some $j \neq i$. Lemma \ref{c(H_i)/c(H_j) < 2} implies that $c_1 \in c(V(G_j))$ or $c_2 \in c(V(G_j))$. Without loss of generality, suppose $c_1 \in c(V(G_j))$ implying $c(V(G_j)) = \{c_1,c_3\}$ by Lemma \ref{|c(Hi)|<3}. It will be shown that $c_1$ appears in every copy of $G$.

\begin{sloppy}Now, for $k \notin \{i,j\}$, Lemma \ref{c(H_i)/c(H_j) < 2} implies that $|c(V(G_i))\setminus c(V(G_k))| \leq 1$ and\end{sloppy}\\ $|c(V(G_j))\setminus c(V(G_k))| \leq 1$. Thus, for all $k \notin \{i,j\}$, either $c_1 \in c(V(G_k))$ or $c_2,c_3 \in c(V(G_k))$ implying $c(V(G_k)) = \{c_2,c_3\}$ by Lemma \ref{|c(Hi)|<3}. Now, define $c': V(H) \to \{red,blue\}$ by 

\[c'(w_k) = 
\begin{cases}
    red & \text{if $c_1 \in c(V(G_k))$,} \\
    blue & \text{if $c(V(G_k)) = \{c_2,c_3\}$.}
\end{cases}\]

For the sake of contradiction, assume $blue \in c'(V(H))$. Then there must exist $red$ and $blue$ neighbors in $H$, call them $w_{\ell_1},w_{\ell_2}$. Without loss of generality, say $c'(w_{\ell_1}) = red$ and $c'(w_{\ell_2}) = blue$ so that $c_1 \in c(V(G_{\ell_1}))$ and $c(V(G_{\ell_2})) = \{c_2,c_3\}$. Then $c_1,c_2,c_3 \in c(V(G_{\ell_1}))\cup c(V(G_{\ell_2}))$ and $3 \leq |c(V(G_{\ell_1}))\cup c(V(G_{\ell_2}))|$, contradicting Corollary \ref{cor:neighborcopies}. Thus, $c'(V(H)) = red$, the desired result.\qed\end{proof}

\section{Graph Products of Paths and Cycles}\label{sec:pmcn}

As a reminder, the conventions for $G\square H$ will be used to label the vertices of $P_m\square C_n$.  In particular, letting $G = P_m$ and $H = C_n$ gives the following:
\begin{itemize}
    \item $V(P_m) = \{u_1,u_2,\dots,u_m\}$ with edges $u_iu_{i+1}$ for $1\le i \le m-1$,
    \item $V(C_n) = \{w_1,w_2,\dots,w_n\}$ with edges $w_iw_{i+1}$ for $1\le i \le n-1$ and $w_nw_1$,
    \item $G_i$ is the $i$th copy of $P_m$ in $P_m\square C_n$ and has vertex set $\{v_{1,i},v_{2,i},\dots ,v_{m,i}\}$, and
    \item $H_i$ is the $i$th copy of $C_n$ in $P_m\square C_n$ and has vertex set $\{v_{i,1},v_{i,2},\dots,v_{i,n}\}$.
\end{itemize}

Now, a fact about $P_m\square C_n$ is that $\dist_{P_m\square C_n}(v_{i,j},v_{k,\ell}) = |i-k| + \min\{ j-\ell \bmod n, \ell - j \bmod n\}$.  Note that the the standard representative of the equivalence class of $\mathbb{Z}_n$ is chosen, i.e. $j-\ell\bmod n, \ell-j \bmod n \in \{0,1,\dots,n-1\}$.

\begin{lemma}\label{lem:p2codd} For any positive integer $k$, $\aw(P_2 \square C_{2k+1},3) = 3$.
\end{lemma}

\begin{proof}
For the sake of contradiction, let $c$ be an exact, rainbow-free $3$-coloring of $P_2 \square C_{2k+1}$. Swapping the roles of $G$ and $H$ in Lemma \ref{|c(Hi)|<3} gives \[|c(V(H_1))|, |c(V(H_2))| \leq 2.\] Without loss of generality, suppose $c(V(H_1)) = \{red,blue\}$, $green \in c(V(H_2))$ with $c(v_{2,1}) = green$.  Note that $c(v_{1,1}) \in \{red,blue\}$ and define $P_\ell$ to be a shortest path in $H_1$ containing $v_{1,1}$ that contains colors $red$ and $blue$.  Without loss of generality, let $P_\ell = (v_{1,1},v_{1,2},v_{1,3},\dots, v_{1,\ell})$ and let $\rho$ be the isometric subgraph of $C_{2k+1}$ that corresponds to $P_\ell$.  Note that $P_2\square\rho$ is an isometric subgraph in $P_2\square C_{2k+1}$ that contains three colors and $\ell \le k+1$.

If $\ell$ is even, then $P_2 \square \rho$ has a rainbow $3$-AP by Theorem \ref{PmxPn} so $P_2\square C_{2k+1}$ has a rainbow $3$-AP, a contradiction.

If $\ell$ is odd and $\ell \le k$, extending $P_\ell$ by one additional vertex (and likewise extending $\rho$ to be $\rho'$)  maintains isometry.  That is, there is an isometric path $P_{\ell+1}$ in $H_1$ that contains $v_{1,1}$ and has colors $red$ and $blue$.  Thus, $P_2\square\rho'$ is an isometric subgraph of $P_2\square C_{2k+1}$ that contains three colors and it contains a rainbow $3$-AP by Theorem \ref{PmxPn}, another contradiction.

Finally, consider the case when $\ell$ is odd and $\ell = k+1$.  Note that $c(v_{1,1}) =red$ for $k+3\le i \le 2k+1$, else the minimality of $P_\ell$ would be contradicted.  Also, $j=\frac{3k+4}{2}$ is an integer and $k+3 \le j \le 2k+1$ when $k\ge 2$.  Thus, $\{v_{2,1},v_{1,j},v_{1,\ell}\}$ is a rainbow $3$-AP, a contradiction.

Therefore, no such $c$ exists and $\aw(P_2 \square C_{2k+1},3) = 3.$\qed\end{proof}

\begin{lemma}\label{lem:pmcodd}
    For integers $m$ and $k$ with $2 \le m$ and $1\le k$, \[\aw(P_m\square C_{2k+1},3) = 3.\]
\end{lemma}

\begin{proof} 
For a base case, note that Lemma \ref{lem:p2codd} implies $\aw(P_2 \square C_{2k+1}, 3) = 3$ for all $1\le k$. As the inductive hypothesis, suppose that $\aw(P_\ell \square C_{2k+1}, 3) = 3$ for some $2\le \ell$. Let $c$ be a rainbow-free, exact $3$-coloring of $P_{\ell+1} \square C_{2k+1}$ and let $H_i$ denote the $i$th copy of $C_{2k+1}$. By hypothesis and the fact that $c$ is rainbow-free, \[\left|c\left(\bigcup_{i=1}^{\ell}V(H_i)\right)\right| \leq 2 \text{ and } \left|c\left(\bigcup_{i=2}^{\ell+1}V(H_i)\right)\right| \leq 2.\]

Thus, the inclusion-exclusion principle gives $\left|c\left(\bigcup_{i=2}^{\ell}V(H_i)\right)\right| = 1$. Without loss of generality, assume \[c\left(\bigcup_{i=2}^{\ell}V(H_i)\right) = \{red\}, \quad blue \in c(V(H_1)), \quad and \quad green \in c\left(V(H_{\ell+1})\right).\]  In particular, assume $c(v_{1,1}) = blue$ and $c(v_{\ell+1,j}) = green$ for some $j \le k+1$.

Suppose $\ell$ is even. Then $\{v_{1,1},v_{\frac{\ell+2}{2},i},v_{\ell+1,j}\}$ is a rainbow $3$-AP for $i = \frac{j+1}{2}$ if $j$ is odd, and $i = \frac{2k+j+2}{2}$ if $j$ is even. On the other hand, suppose $\ell$ is odd.
Then, $\{v_{1,1},v_{\frac{\ell+1}{2},i},v_{\ell+1,j}\}$ is a rainbow $3$-AP for $i = \frac{j+2}{2}$ if $j$ is even, and $i = \frac{2k+j+1}{2}$ if $j$ is odd. 

In any case, there is a rainbow $3$-AP which is a contradiction, so $\aw(P_{\ell+1} \square C_{2k+1}, 3) = 3.$  Thus, by induction, $\aw(P_{m} \square C_{2k+1}, 3) = 3$ for any $2\le m$.\qed\end{proof}

Determining of $\aw(P_2\square C_{2k})$ requires two strategies since there are $k$ values for which $\aw(P_2\square C_{2k})=3$ and $k$ values for which $\aw(P_2\square C_{2k})=4$. Essentially, $\aw(P_2\square C_n)=4$ when $n=4\ell$ and is determined by providing a coloring where one pair of vertices that are diametrically opposed are colored distinctly and everything else is a third color. This avoids rainbow $3$-APs since the diameter of $P_2\square C_{4\ell}$ is odd and because each vertex $v\in V(C_{4\ell})$ has exactly one vertex whose distance from $v$ realizes the diameter of $C_{4\ell}$.  Note that this is different than what happens in $P_2\square C_{2k+1}$ since each vertex $v \in V(C_{2k+1})$ has two vertices whose distance from $v$ realizes the diameter of $C_{2k+1}$. When the diameter of $P_2\square C_{2k}$ is even, this coloring, and every other coloring, ends up creating an isometric $P_2\square P_{2j}$ with $3$-colors.  Then, it is only a matter of applying Theorem \ref{PmxPn} to find the rainbow $3$-AP.

\begin{lemma}\label{lem:pmcevendiamodd}
     \begin{sloppy}
     For integers $m$ and $k$ with $2 \le m,k$, $\aw(P_m \square C_{2k},3) = 4$ if $\diam(P_m \square C_{2k})$ is odd.
     \end{sloppy}
\end{lemma}

\begin{proof}
    Define $c: V(P_m \square C_{2k}) \to \{red,blue,green\}$ by 
    \[c(v_{i,j}) = 
    \begin{cases}
        blue & \text{if $i=j=1$,} \\
        green & \text{if $i=m$ and $j=k+1$,} \\
        red & \text{otherwise.}
    \end{cases}\]
    Note that any rainbow $3$-AP must contain $v_{1,1}$ and $v_{m,k+1}$ since they are the only $blue$ and $green$ vertices, respectively. This will be shown by proving $v_{1,1}$ and $v_{m,k+1}$ are not part of any nondegenerate $3$-AP. For the sake of contradiction, assume there exists $v_{i,j}\in V(P_m\square C_n)$ such that $\{v_{1,1},v_{i,j},v_{m,k+1}\}$ is a nondegenerate $3$-AP. 
    
    One way this can happen is if $\dist(v_{1,1},v_{i,j}) = \dist(v_{i,j},v_{m,k+1})$. Without loss of generality, suppose $1 \leq j \leq k+1$. Then \[\dist(v_{1,1},v_{i,j}) = (i-1) + (j-1) = i+j-2,\] and \[\dist(v_{i,j},v_{m,k+1}) = (m-i) + (k+1-j) = m+k+1-i-j.\] By assumption, $i+j-2 = m+k+1-i-j$ which implies that $m+k-1 = 2i+2j-2$. However, $\diam(P_m \square C_{2k}) = m+k-1$ is odd, a contradiction.  
    
    The only other possible way that $\{v_{1,1},v_{i,j}, v_{m,k+1}\}$ is a $3$-AP is if $\dist(v_{i,j},v_{1,1}) = \diam(P_m\square C_{2k})$ or $\dist(v_{i,j},v_{m,k+1}) = \diam(P_m\square C_{2k})$.  However, this implies $v_{i,j}\in\{v_{1,1},v_{m,k+1}\}$ which gives a degenerate $3$-AP.
    
    Thus, the exact $3$-coloring $c$ of $P_m\square C_{2k}$ is rainbow free so $4\le \aw(P_m \square C_{2k},3)$.  Theorem \ref{thm:rsw} gives an upper bound of $4$ which implies $\aw(P_m \square C_{2k},3) = 4$.\qed\end{proof}

\begin{lemma}\label{lem:p2ceven}
For any integer $k$ with $2\le k$, 
\[\aw(P_2 \square C_{2k},3) = \begin{cases}
3 & \text{if $k$ is odd,} \\
4 & \text{if $k$ is even.}
\end{cases}\]
\end{lemma}

\begin{proof}
If $k$ is even, then $\diam(P_2\square C_{2k}) = 1+k$ is odd so by Lemma \ref{lem:pmcevendiamodd} $\aw(P_2\square C_{2k}) = 4$.

Now assume $k$ is odd and let $c$ be an exact $3$-coloring of $P_2 \square C_{2k}$. For the sake of contradiction, assume $c$ is rainbow-free. By Lemma \ref{|c(Hi)|<3}, $|c(V(H_1))|,|c(V(H_2))| \leq 2$. Without loss of generality, suppose $c(V(H_1)) = \{red,blue\}$, and $green \in c(V(H_2))$ with $c(v_{2,1}) = green$.  Now, define $P_\ell$ as a shortest path in $H_1$ containing $v_{1,1}$ that contains colors $red$ and $blue$, and let $\rho$ be the isometric subgraph of $C_{2k+1}$ that corresponds to $P_\ell$. Note that $P_2\square \rho$ is an isometric subgraph in $P_2\square C_{2k}$ that contains three colors. 
If $\ell$ is even, then Theorem \ref{PmxPn} gives a rainbow $3$-AP, a contradiction.  

If $\ell$ is odd, then extending $P_\ell$ by one vertex in either direction maintains isometry.  In other words, there is an isometric path $P_{\ell+1}$ in $H_1$ that contains $v_{1,j}$ and the colors $red$ and $blue$.  Thus, $P\square P_{\ell+1}$ is an isometric subgraph of $P_2\square C_{2k}$ that contains three colors which means it has a rainbow $3$-AP by Theorem \ref{PmxPn}, a contradiction.

Therefore, when $k$ is odd, every exact $3$-coloring of $P_2\square C_{2k}$ has a rainbow $3$-AP and $\aw(P_2 \square C_{2k},3) = 3$.\qed\end{proof}

Before getting to more general results an analysis of $\aw(P_3\square C_n)$ needs to happen.  Similar to the $\aw(P_2\square C_n)$ situation, there are very subtle and important differences when $n$ is odd versus when $n$ is even.

\begin{lemma}\label{lem:p3ceven}
    For any integer $k$ with $2\le k$, \[\aw(P_3\square C_{2k},3) = \begin{cases}
    3 & \text{if $k$ is even,} \\
    4 & \text{if $k$ is odd.}
    \end{cases}\]
\end{lemma}

\begin{proof} 
    If $k$ is odd, then $\diam(P_3\square C_{2k}) = 2+k$ is odd so by Lemma \ref{lem:pmcevendiamodd} $\aw(P_2\square C_{2k}) = 4$.

    Suppose $k$ is even and $c$ is an exact, rainbow-free $3$-coloring of $P_3\square C_{2k}$. Then an argument similar to the argument in the proof of Lemma \ref{lem:pmcodd} can be used to establish, without loss of generality, that $c(V(H_1)) = \{red,blue\}$, $c(V(H_2)) = \{red\}$, $c(V(H_3)) = \{red,green\}$, $c(v_{1,1}) = blue$ and $c(v_{3,j}) = green$ for some $1\le j \le k+1$. 
  
    If $j$ is odd, then $\{v_{1,1},v_{2,\frac{j+1}{2}},v_{3,j}\}$ is a rainbow $3$-AP, contradicting that $P_3\square C_{2k}$ is rainbow free. So, suppose $j$ is even. Then $j+1\leq k+1$ implying that the path $P_{j+1} = (w_1,\ldots,w_{j+1})$ is an isometric subgraph of $C_{2k}$. So, $P_3\square P_{j+1}$ is an isometric subgraph of $P_3\square C_{2k}$. Since $c(P_3\square P_{j+1}) = \{red,blue,green\}$, Theorem \ref{PmxPn} implies that $P_3\square C_{2k}$ contains a rainbow $3$-AP.\qed\end{proof}

\begin{lemma}\label{lem:pmc4x+2}
    If $m \geq 2$ is even and $k\geq 1$, then \[\aw(P_m\square C_{4k+2},3) = 3.\]
\end{lemma}

\begin{proof}
Lemma \ref{lem:p2ceven} implies $\aw(P_2 \square C_{4k+2},3) = 3$. Suppose $\aw(P_\ell \square C_{4k+2}, 3) = 3$ for some even $\ell \geq 2$. Then, let $c$ be an exact 3-coloring of $P_{\ell+2} \square C_{4k+2}$ that avoids rainbow $3$-APs, and let $H_i$ denote the $i$th copy of $C_{4k+2}$.  By hypothesis, \[\left|c\left(\bigcup_{i=1}^{\ell}V(H_i)\right)\right| \leq 2 \quad \text{and} \quad  \left|c\left(\bigcup_{i=3}^{\ell+2}V(H_i)\right)\right| \leq 2.\] By the inclusion-exclusion principle, $\left|c\left(\bigcup_{i=3}^{\ell}V(H_i)\right)\right| = 1$. Without loss of generality, suppose $c\left(\bigcup_{i=3}^{\ell}V(H_i)\right) = \{red\}$, so that Proposition \ref{prop:everycopy} implies $red\in c(H_i)$ for $1\le i \le \ell+2$. Further, without loss of generality, suppose $blue \in c(V(H_1)\cup V(H_2))$ and $green \in c(V(H_{\ell+1})\cup V(H_{\ell+2}))$. Say, $c(v_{i,1}) = blue$ and $c(v_{h,j}) = green$ for $i \in \{1,2\}$, $h \in \{\ell+1,\ell+2\}$ and $1 \leq j \leq 2k+1$ such that $i$ is maximal and $h$ is minimal.  If $i=2$ and $h=3$, then $|c(H_2)\cup c(H_3)| \ge 3$ which contradictions Lemma \ref{cor:neighborcopies}.  So assume $h-i \ge 2$. Thus, $c(V(H_{i+1})) = \{red\}$ and $c(V(H_{h-1})) = \{red\}$. 

\begin{description}
\item[Case 1.] Suppose $\dist(v_{i,1},v_{h,j})$ is even. 

Then either $\dist_{P_{\ell+2}}(u_i,u_h)=h-i$ and $\dist_{C_{4x+2}}(w_1,w_j)=j-1$ are both odd or both even. If they are both even, then $\{v_{i,1}, v_{\frac{i+h}{2},\frac{j+1}{2}}, v_{h,j}\}$ is a rainbow $3$-AP. If they are both odd, then $\{v_{i,1}, v_{\frac{i+h+1}{2},\frac{j}{2}},v_{h,j}\}$ is a rainbow $3$-AP.
\item[Case 2.] Suppose $\dist(v_{i,1},v_{h,j})$ is odd. 

If  $j < 2k+2$, then $\{v_{h,j},v_{i,1},v_{h-1,j+1}\}$ is a rainbow $3$-AP. So, suppose $j=2k+2$. Then $\dist_{C_{4k+2}}(w_1,w_j)=2k+1$ is odd implying that $\dist_{P_{\ell+2}}(u_i,u_h)$ is even. Thus, either $i = 1$ and $h = \ell+1$, or $i = 2$ and $h = \ell + 2$. First, suppose $i = 1$ and $h = \ell + 1$. Then the $3$-AP $\{v_{\ell+1,j},v_{1,1},v_{\ell+2,j+1}\}$ implies $c(v_{\ell+2,j+1}) = green$. Since $i$ is maximal, $c(V(H_2)) = \{red\}$. Thus, $\{v_{1,1}, v_{\ell+2,j+1}, v_{2,2}\}$ is a rainbow $3$-AP since $j+1 = 2k+3$. For $i=2$ and $j=\ell+2$, the $3$-APs $\{v_{2,1},v_{\ell+2,j},v_{1,2}\}$ and $\{v_{\ell+2,j},v_{1,2},v_{\ell+1,j+1}\}$ yield a rainbow $3$-AP.
\end{description}
Thus, $\aw(P_{\ell+2} \square C_{4k+2},3)=3$ and by induction, $\aw(P_{m} \square C_{4k+2},3)=3$ for any even $m\geq2$.\qed\end{proof}

Replacing $4k+2$ with $4k$ and $2k+2$ with $2k+1$ gives the proof of Lemma \ref{lem:pmc4x}, thus the proof has been omitted.

\begin{lemma}\label{lem:pmc4x}
If $m \geq 3$ is odd and $k\geq 1$, then \[\aw(P_m\square C_{4k},3) = 3.\]
\end{lemma}

Lemmas \ref{lem:pmcodd}, \ref{lem:pmcevendiamodd}, \ref{lem:pmc4x+2}, and \ref{lem:pmc4x} yield the following theorem.

\begin{theorem}\label{thm:pmcn}
    If $m\geq2, n\geq3$ then \[\aw(P_m\square C_n,3)= 
    \begin{cases}
        4 & \text{if $n$ is even and $\diam(P_m\square C_n)$ is odd,} \\
        3 & \text{otherwise.}
    \end{cases}
    \]
\end{theorem}

\section{Graph Products of Cycles with Other Graphs}\label{sec:cmcn}

This section starts with a general result, Theorem \ref{CnxG}, and then uses the general result to establish $\aw(C_m\square C_n,3)$.

\begin{theorem}\label{CnxG}
For any integer $k$ with $1\le k$, $\aw(G \square C_{2k+1}, 3) = 3$ for any connected graph $G$ with $|G|\geq 2$. 
\end{theorem}

\begin{proof}
    Let $V(G) = \{u_1,\ldots,u_n\}$ and $H_i$ denote the $i$th labeled copy of $C_{2k+1}$.
    Lemma \ref{lem:p2codd} implies that $\aw(P_2\square C_{2k+1},3)=3$, so suppose $|G|\ge3$.  Let $c:V(G \square C_{2k+1}) \to \{red,blue,green\}$ be an exact $3$-coloring, and, for the sake of contradiction, assume $c$ is rainbow-free. Since $|G|\geq 3$, Proposition \ref{prop:everycopy} implies that, without loss of generality, $red$ is in every copy of $C_{2k+1}$. So, define $c': V(G) \to \{red,blue,green\}$ by 
    \[c'(u_i) = 
    \begin{cases}
       red & \text{if $c(V(H_i)) = \{red\}$,} \\
        \C & \text{if $\C \in c(V(H_i))\setminus\{red\}$.}
    \end{cases}\]
    Since Lemma \ref{|c(Hi)|<3} implies that $|c(V(H_i))| \leq 2$ for all $1\leq i\leq n$, it follows that $c'$ is well-defined. By Lemma \ref{isometricpathorC3}, there either exists a $C_3$ in $G$ containing $red$, $blue$, and $green$ or an isometric path in $G$ containing $red$, $blue$, and $green$. 
    \par
    First, suppose $C_3 \cong G[\{u_{i_1},u_{i_2},u_{i_3}\}]$ contains $red$, $blue$, and $green$. Then, without loss of generality, there exists neighboring copies $H_{i_1}$ and $H_{i_2}$ of $H$, in $G\square C_{2k+1}$, such that $c(V(H_{i_1})) = \{red,blue\}$ and $c(V(H_{i_2})) = \{red,green\}$, contradicting Corollary \ref{cor:neighborcopies}.
    \par
    Finally, suppose there exists an isometric path $P$ in $G$ such that $c'(V(P)) = \{red,blue,\\green\}$. Now, by Lemma \ref{lem:pmcodd}, there exists a rainbow $3$-AP in the isometric subgraph $P \square C_{2k+1}$, a contradiction.\qed\end{proof}

Just as Lemma \ref{lem:pmcodd} was generalized into Theorem \ref{CnxG} which showed that \[\aw(G \square C_{2k+1},3)=3\] for all connected $G$ with at least $2$ vertices, significant time was spent on the conjecture that a similar generalization could be performed to show $\aw(G \square C_{4k+2},3)=3$ when $\diam(G)$ is odd and $\aw(G \square C_{4k},3)=3$ when $\diam(G)$ is even. However, these conjectures do not hold because it cannot be guaranteed that an isometric $P_{2j} \square C_{4k+2}$ subgraph of $G \square C_{4k+2}$ or $P_{2j+1} \square C_{4k}$ subgraph of $G \square C_{4k}$ existed that contained three colors. The following example provides such a $G$.

\begin{example}\label{ex:counterex}
Consider the graph in Figure \ref{fig:pncevencex} which is $G \square C_4$, where $G$ is a $C_{10}$ with a leaf. That is $V(G) = \{w_1,\ldots,w_{11}\}$ with edges $w_iw_{i+1}$ for $1\le i \le 9$ and the additional edges $w_1w_{10}$ and $w_{10}w_{11}$. Define $c:V(G\square C_4) \to \{red,blue,green\}$ by $c(v_{2,1}) = blue$, $c(v_{7,3}) = green$, and $c(v) = red$ for all $v \in V(G \square C_4)\setminus\{v_{2,1},v_{7,3}\}$. In order for $G \square C_4$ to contain a rainbow $3$-AP, there must exist a red $v \in V(G \square C_4)$ such that \[\dist(v_{2,1},v) = \dist(v,v_{7,3}), \quad \dist(v,v_{2,1}) = \dist(v_{2,1},v_{7,3}), \quad \text{or} \quad \dist(v,v_{7,3}) = \dist(v_{7,3},v_{2,1}).\] By construction, every vertex $v$ of $G \square C_4$ is such that $\dist(v,v_{2,1})$ and $\dist(v,v_{7,3})$ have different parity, thus $\dist(v_{2,1},v) \neq \dist(v,v_{7,3})$ for all $v \in V(G)$.
To show that there are no vertices $v$ of $G$ distinct from $v_{2,1},v_{7,3}$ such that $\dist(v,v_{2,1}) = \dist(v_{2,1},v_{7,3})$ or $\dist(v,v_{7,3}) = \dist(v_{7,3},v_{2,1})$, a discussion about \emph{eccentricity} is needed.  For a vertex $v$ of a graph $G$, the \textit{eccentricity} of $v$, denoted $\epsilon(v)$, is the distance between $v$ and a vertex furthest from $v$ in $G$. In other words, \[\epsilon(v) = \max_{u\in V(G)} \dist(u,v).\]  In this example, $\epsilon(v_{2,1}) = \epsilon(v_{7,3}) = \dist(v_{2,1},v_{7,3}) = 7$ and both eccentricities are unique-ly realized. So, there are no non-degenerate $3$-APs in $G \square C_4$ containing $v_{2,1}$ and $v_{7,3}$. Thus, $\aw(G \square C_4,3)=4$. 

\begin{figure}[ht!]
    \centering
    \def \scale {.5} 
    \begin{tikzpicture}
        \def \n {10}
    	\def \radius {5*\scale}
    	\def \radiuss {7*\scale}
    	\def \radiusss {6.1*\scale}
    	
    	\foreach \s in {1,...,\n}
    	{					
    	    \node[draw, circle, fill = red, opacity = .6, inner sep = 3.5*(\scale)^.5] (1\s) at (({360/\n * (\s-1)) + 90}:\radius) {};
    	    \node[draw, circle, fill = red, opacity = .6, inner sep = 3.5*(\scale)^.5] (3\s) at (({360/\n * (\s-1)) + 90}:\radiuss) {};
    	    \node[draw, circle, fill = red, opacity = .6, inner sep = 3.5*(\scale)^.5] (2\s) at (({360/\n * (\s-1)) + 99.25}:\radiusss) {};
    	    \node[draw, circle, fill = red, opacity = .6, inner sep = 3.5*(\scale)^.5] (4\s) at (({360/\n * (\s-1)) + 80.75}:\radiusss) {};
        }
        
        \node[draw, circle, fill = red, opacity = .6, inner sep = 3.5*(\scale)^.5] (111) at (-.7*\scale,3.5*\scale) {};
        \node[draw, circle, fill = red, opacity = .6, inner sep = 3.5*(\scale)^.5] (211) at (-1.7*\scale,2.5*\scale) {};
        \node[draw, circle, fill = red, opacity = .6, inner sep = 3.5*(\scale)^.5] (311) at (-.7*\scale,1.5*\scale) {};
        \node[draw, circle, fill = red, opacity = .6, inner sep = 3.5*(\scale)^.5] (411) at (.3*\scale,2.5*\scale) {};
        
        \foreach \s in {1,...,11}
        {
            \draw (1\s) to (2\s);
    	    \draw (2\s) to (3\s);
    	    \draw (3\s) to (4\s);
    	    \draw (4\s) to (1\s);
        }
        
        \foreach \r in {1,...,4}
        {
            \draw (\r1) to (\r2);
            \draw (\r2) to (\r3);
            \draw (\r3) to (\r4);
            \draw (\r4) to (\r5);
            \draw (\r5) to (\r6);
            \draw (\r6) to (\r7);
            \draw (\r7) to (\r8);
            \draw (\r8) to (\r9);
            \draw (\r9) to (\r10);
            \draw (\r10) to (\r1);
            \draw (\r1) to (\r11);
        }
        
        \node[draw, circle, fill = white, opacity = 1, inner sep = 3.5*(\scale)^.5] () at (({360/10 * (2)) + 90}:\radiuss) {};
        
        \node[draw, circle, fill = white, opacity = 1, inner sep = 3.5*(\scale)^.5] () at (({360/10 * (7)) + 90}:\radius) {};
        
        \node[draw, circle, fill = blue, opacity = .5, inner sep = 3.5*(\scale)^.5] () at (({360/10 * (2)) + 90}:\radiuss) {};
        
        \node[draw, circle, fill = green, opacity = .7, inner sep = 3.5*(\scale)^.5] () at (({360/10 * (7)) + 90}:\radius) {};
        
        \node at (({360/10 * (2)) + 91.5}:8*\scale) {$v_{2,1}$};
        
        \node at (({360/10 * (7)) + 88.5}:4*\scale) {$v_{7,3}$};
    \end{tikzpicture}
    \caption{Image for Example \ref{ex:counterex}: Graph $G\square C_4$, counterexample of generalizing Lemma \ref{lem:pmc4x}.}
    \label{fig:pncevencex}
\end{figure}
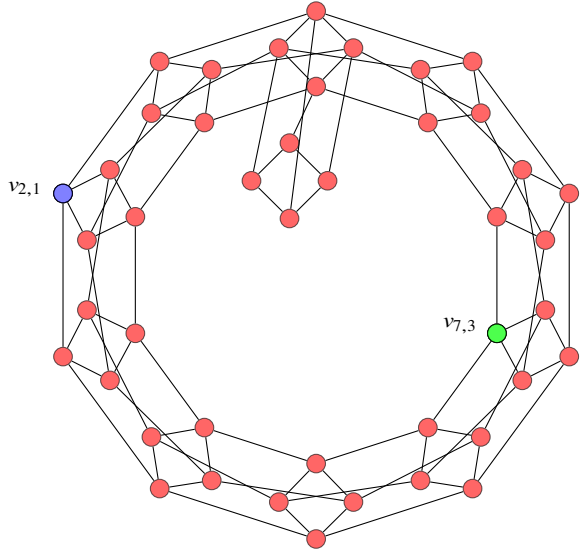
\end{example}

Note that the graph in Figure \ref{fig:pncevencex} is the only example presented in this paper of a graph product with even diameter and anti-van der Waerden number (with respect to $3$) equal to $4$.  This is discussed more in Section \ref{sec:future}.

Theorem \ref{CnxG} gives the following result.

\begin{corollary}\label{cor:coddcodd}
    If $m$ or $n$ is odd with $m,n \geq 3$, then $\aw(C_m\square C_n,3) = 3$.
\end{corollary}

Lemmas \ref{lem:pmc4x+2} and \ref{lem:pmc4x} are used to prove Lemma \ref{lem:cevencevenaw3}.

\begin{lemma}\label{lem:cevencevenaw3}
    If $m$ and $n$ are even with $m\equiv n \Mod{4}$, then $\aw(C_m\square C_n,3) = 3$. 
\end{lemma}

\begin{proof}
    Let $c$ be an exact $3$-coloring of $C_m\square C_n$. Lemma \ref{isometricpathorC3} implies that $C_m\square C_n$ either contains an isometric path or a $C_3$ with three colors. Since there are no $C_3$ subgraphs in $C_m\square C_n$, it follows that $C_m\square C_n$ must contain an isometric path with three colors. Call a shortest such path $P$. Suppose $P$ intersects $k$ copies of $C_n$, and, without loss of generality, suppose these copies are $H_1,\ldots,H_k$. 
    \par
    Notice that there are vertices $v$ and $v'$ of $P$ in $V(H_1)$ and $V(H_k)$, respectively. If $k > \frac{m}{2} + 1$, then any shortest path from $v$ to $v'$ would be contained in the subgraph induced by the vertices of $H_k,H_{k+1},\ldots,H_n,H_1$. So, no shortest path between $v$ and $v'$ would be contained in $P$, implying that $P$ is not isometric, a contradiction. 
    \par
    Thus, $k \leq \frac{m}{2} + 1$, and $P$ is a subgraph of $P_{\frac{m}{2}+1}\square C_n$ where $P_{\frac{m}{2}+1}$ is the subgraph of $C_m$ induced by $\{u_1,\ldots,u_{\frac{m}{2}+1}\}$. Thus, $P$ is an isometric subgraph of $C_m\square C_n$ because $P_{\frac{m}{2} + 1}$ is isometric in $C_m$. Since there are three colors in $P$, there are three colors in $P_{\frac{m}{2} + 1} \square C_n$. Furthermore, since $m\equiv n \Mod{4}$, $\frac{m}{2}$ and $\frac{n}{2}+1$ have different parity. So, Lemma \ref{lem:pmc4x+2} or Lemma \ref{lem:pmc4x} implies that $P_{\frac{m}{2} + 1} \square C_n$ contains a rainbow $3$-AP. Thus, $C_m\square C_n$ contains a rainbow $3$-AP.\qed\end{proof}

In the proof of Lemma \ref{lem:cevencevenaw4}, the fact that each vertex in an even cycle realizes the diameter with exactly one other vertex will be used.

\begin{lemma}\label{lem:cevencevenaw4}
    If $m$ and $n$ are even with $m\not\equiv n \Mod{4}$, then $\aw(C_m\square C_n,3) = 4$.
\end{lemma}

\begin{proof}
    Define $k=\frac{m}{2}+1$ and $\ell=\frac{n}{2}+1$ and the coloring $c: V(C_m\square C_n)\to \{red,blue,\\green\}$ by
    \[c(v_{i,j}) = 
    \begin{cases}
        blue & \text{if $i=j=1$,} \\
        green & \text{if $i=k, j=\ell$,} \\
        red & \text{otherwise.}
    \end{cases}\]
    Since $v_{1,1}$ and $v_{k,\ell}$ are the only $blue$ and $green$ vertices, any rainbow $3$-AP must contain them. This result will be proved by showing $v_{1,1}$ and $v_{k,\ell}$ are not part of any nondegenerate $3$-AP. For the sake of contradiction, assume there exists $v_{i,j}\in V(C_m\square C_n)$ such that $\{v_{1,1},v_{i,j},v_{k,\ell}\}$ is a nondegenerate $3$-AP. 
    
    One way this can happen is if $\dist(v_{1,1},v_{i,j}) = \dist(v_{i,j},v_{k,\ell})$. Without loss of generality, up to a relabelling of the vertices, suppose $1 \leq i \leq k$ and $1 \leq j \leq \ell$. Then \[\dist(v_{1,1},v_{i,j}) = (i-1) + (j-1) = i+j-2,\] and \[\dist(v_{i,j},v_{k,\ell}) = (k-i) + (\ell-j) = k+\ell-i-j.\] By assumption, $i+j-2 = k+\ell-i-j$, which implies that 
    \begin{equation}\label{eq1}
        2i + 2j - 2 = k+\ell = \frac{m}{2} + \frac{n}{2}.
    \end{equation}
    However, $m\not\equiv n \Mod{4}$ implies $\frac{m}{2} + \frac{n}{2}$ is odd, which contradicts equation (\ref{eq1}).
    
    The only other possible way that $\{v_{1,1},v_{i,j}, v_{k,\ell}\}$ is a $3$-AP is if $\dist(v_{i,j}, v_{1,1}) = \dist(v_{1,1},v_{k,\ell})$ or $\dist(v_{i,j}, v_{k,\ell}) = \dist(v_{k,\ell},v_{1,1})$.  However, \[\epsilon(v_{1,1})=\epsilon(v_{k,\ell})=\diam(C_m\square C_n)\] is uniquely realized. This implies $v_{i,j}\in \{v_{1,1},v_{k,\ell}\}$ yielding a degenerate $3$-AP.
    
    Thus, the exact $3$-coloring $c$ of $C_m\square C_n$ is rainbow free so $4\le \aw(C_m\square C_n,3)$.  Theorem \ref{thm:rsw} gives an upper bound of $4$ which implies $\aw(C_m\square C_n,3) = 4$.
\qed\end{proof}

Conglomerating Corollary \ref{cor:coddcodd}, Lemma \ref{lem:cevencevenaw3} and Lemma \ref{lem:cevencevenaw4} yields Theorem \ref{thm:cmcn}.

\begin{theorem}\label{thm:cmcn}
    If $m,n \geq 3$, then
    \[\aw(C_m\square C_n,3) = 
    \begin{cases}
        4 & \text{if $m$ and $n$ are even and $\diam(C_m\square C_n)$ is odd,} \\
        3 & \text{otherwise.}
    \end{cases}
    \]
\end{theorem}

\section{Future Work}\label{sec:future}

Recall that Example \ref{ex:counterex} was the only example  presented in this paper of a graph product with even diameter and anti-van der Waerden number (with respect to $3$) equal to $4$. One of the key factors in allowing this to happen was a pair of vertices $u$ and $v$ such that $\epsilon(u)=\epsilon(v)=\dist(u,v)<\diam(u,v)$. Such vertices will be called \textit{almost peripheral vertices} whose name comes from \textit{peripheral vertices} which are vertices who realize the diameter.

\begin{conjecture}\label{conj:5.1}
    If $G \square H$ has no almost peripheral vertices and $\diam(G\square H)$ is even, then $\aw(G\square H,3)=3$.
\end{conjecture}
In particular, the authors believe that trees do not contain any almost peripheral vertices. For this reason, it is believed that Conjecture \ref{conj:5.2} holds if Conjecture \ref{conj:5.1} holds.
\begin{conjecture}\label{conj:5.2}
    If $T$ is a tree, $n$ is even, and $\diam(T\square C_n)$ is even, then $\aw(T\square C_n,3)=3$. 
\end{conjecture}
This result would provide a more specific case of when the even cycle analog of Theorem \ref{CnxG} holds.

Another way to extend Theorem \ref{CnxG} would be considering $\aw(G \square C_n,k)$ for some $k>3$. For, $k=3$, Theorem \ref{CnxG} showed that when $n$ is odd, $\aw(G \square C_n,k)=k$ for any connected $G$ of order at least $2$. However, there may be other properties of $n$ that guarantee $\aw(G \square C_n,k)=k$ for $k > 3$. Some preliminary work analyzing $\aw(P_m\square C_n,4)$ suggests that for any $n$, there exists an $m$ such that $\aw(P_m\square C_n,4)\geq 5$.

\begin{acknowledgements}
Thanks to Ethan Manhart, Hunter Rehm and Laura Zinnel for providing feedback and offering ideas during this project.
\end{acknowledgements}


\begin{thebibliography}{20}

    \bibitem{AF} M. Axenovich and D. Fon-Der-Flaass, On rainbow arithmetic progressions, \emph{Electron. J. Combin.} {\bf 11} (2004), no. 1, Research Paper 1, 7pp.
    
    \bibitem{AM} M. Axenovich and R.R. Martin. Sub-Ramsey numbers for arithmetic progressions. \emph{Graphs and Combinatorics,} {\bf 22} (2006), no. 1, 297--309.
		
	\bibitem{SWY} Z. Berikkyzy, A. Schulte, E. Sprangel, S. Walker, N. Warnberg and M. Young, Anti-van der Waerden numbers on Graphs, Accepted to \emph{Graphs and Combinatorics}, \url{https://arxiv.org/pdf/1802.01509.pdf}.

    \bibitem{BSY} Z. Berikkyzy, A. Schulte, and M. Young, Anti-van der Waerden numbers of 3-term arithmetic progressions, \emph{Electron. J. Combin. } {\bf 24} (2017), no. 2, Paper 2.39, 9 pp.
	
	\bibitem{BKKTTY} E. Bevilacqua, A. King, J. Kritschgau, M. Tait, S. Tebon and M. Young, Rainbow numbers for $x_1 + x_2 = kx_3$ in $\mathbb{Z}_n$, \emph{Integers} {\bf 20} (2020),  A50.
	
    \bibitem{DMS} S. Butler, C. Erickson, L. Hogben, K. Hogenson, L. Kramer, R.L. Kramer, J. Lin, R.R. Martin, D. Stolee, N. Warnberg and M. Young, Rainbow Arithmetic Progressions, \emph{J. Comb.} {\bf 7} (2016), no. 4, 595--626.
    
    \bibitem{ESS} P. Erd{\H{o}}s, M. Simonovits, and V. S{\'{o}}s Anti-Ramsey theorems. \emph{Infinite and finite sets (Colloq., Keszthely, 1973; dedicated to P. Erd{\H{o}}s on his 60th birthday)} {\bf{II}} (1973), 633–643.

    \bibitem{FGRWW} K. Fallon, C. Giles, H. Rehm, S. Wagner and N. Warnberg, Rainbow numbers of $[n]$ for $\sum_{i=1}^{k-1} x_i = x_k$, \emph{Austral. J. Combin.} {\bf 77(1)} (2020), 1--8.
    
    \bibitem{FMO} S. Fujita, C. Magnant, and K. Ozeki, Rainbow generalizations of Ramsey theory: A dynamic survey, \emph{Theory Appl. Graphs,} {\bf 0} (2014), no. 1, Article 1.


    \bibitem{RFC} M. Huicochea and A. Montejano, The Structure of Rainbow-Free Colorings For Linear Equations on Three Variables in $\mathbb{Z}_p$, \emph{Integers} {\bf 15A} (2015), A8.

    \bibitem{J} V. Jungi\'c, J. Licht (Fox), M. Mahdian, J. Ne\u{s}etril, and R. Radoi\u{c}i\'c, Rainbow arithmetic progressions and anti-Ramsey results, \emph{Combinatorics, Probability and Computing} {\bf 12} (2003), no 5-6, 599--620.
	
	\bibitem{LM} B. Llano and A. Montejano, Rainbow-free Colorings of $x+y = cz$ in $\mathbb{Z}_p$, \emph{Discrete Math.} {\bf 312}, (2012), 2566--2573.
	
	\bibitem{R} F.P. Ramsey,  On a Problem of Formal Logic, \emph{Proc. London Math. Soc.} {\bf 30} (1928), 264 - 286.
	
	\bibitem{RSW} H. Rehm, A. Schulte and N. Warnberg, Anti-van der Waerden numbers on Graph Products, \emph{Austral. J. Combin.} {\bf 73(3)} (2019), 486--500.	
	
	\bibitem{S} I. Schur, \"{U}ber Potenzreihen die im Innern des Einheitskreises beschr\"{a}nkt sind. J. Reine Angew. Math (1917), 205-232.
		
    \bibitem{U} K. Uherka. An introduction to Ramsey theory and anti-Ramsey theory on the integers. Master's Creative Component (2013), Iowa State University.
		
    \bibitem{W27} B. van der Waerden, Beweis einer baudetschen vermutung. {\emph{Nieuw Arch. Wisk.}} {\bf{19}} (1927), 212–216.

    \bibitem{finabgroup} M. Young.  Rainbow Arithmetic Progressions in Finite Abelian Groups, \emph{J. Comb.} {\bf 9(4)} (2018), 619--629.
\end{thebibliography}
\end{document}